\documentclass[journal,twoside,web]{ieeecolor}
\usepackage{lcsys}
\usepackage{cite}
\usepackage{amsmath,amssymb,amsfonts,mathrsfs, mathtools}
\usepackage{algorithmic, tikz}
\usepackage{graphicx}
\usepackage{textcomp}
\newtheorem{theorem}{Theorem}
\newtheorem{lemma}{Lemma}
\newtheorem{assumption}{Assumption}
\newtheorem{remark}{Remark}

\newtheorem{definition}{Definition}
\def\BibTeX{{\rm B\kern-.05em{\sc i\kern-.025em b}\kern-.08em
    T\kern-.1667em\lower.7ex\hbox{E}\kern-.125emX}}
\newcommand\copyrighttext{%
	\footnotesize \copyright 2023 IEEE. Personal use of this material is permitted. Permission from IEEE must be obtained for all other uses, in any current or future media, including reprinting/republishing this material for advertising or promotional purposes, creating new collective works, for resale or redistribution to servers or lists, or reuse of any copyrighted component of this work in other works.}
\newcommand\copyrightnotice{%
	\begin{tikzpicture}[remember picture,overlay]
		\node[anchor=south,yshift=5pt] at (current page.south) {\fbox{\parbox{\dimexpr\textwidth-\fboxsep-\fboxrule\relax}{\copyrighttext}}};
	\end{tikzpicture}%
}

\begin{document}
\title{\LARGE On the design of persistently exciting inputs for data-driven control of linear and nonlinear systems}

\author{Mohammad Alsalti, Victor G. Lopez, and Matthias A. Müller %
	\thanks{Leibniz University Hannover, Institute of Automatic Control, 30167 Hannover, Germany. E-mail:\{alsalti,lopez,mueller\}@irt.uni-hannover.de.}%
	\thanks{This work has received funding from the European Research Council (ERC) under the European Union’s Horizon 2020 research and innovation programme (grant agreement No 948679).
	}
}

\maketitle
\thispagestyle{empty}
\copyrightnotice
\begin{abstract}%
	In the context of data-driven control, persistence of excitation (PE) of an input sequence is defined in terms of a rank condition on the Hankel matrix of the input data. For nonlinear systems, recent results employed rank conditions involving collected input and state/output data, for which no guidelines are available on how to satisfy them a priori. In this paper, we first show that a set of discrete impulses is guaranteed to be persistently exciting for any controllable LTI system. Based on this result, for certain classes of nonlinear systems, we guarantee persistence of excitation of sequences of basis functions \textit{a priori}, by design of the physical input only.
\end{abstract}


	\section{Introduction}
\IEEEPARstart{C}{}entral to the area of data-driven control, as well as system identification \cite{Ljung87} and adaptive control \cite{astrom08}, is the notion of persistence of excitation (PE), see, e.g., \cite{Green86}. In \cite{Willems05}, a discrete-time sequence is said to be persistently exciting of a certain order if a corresponding Hankel matrix of that sequence has full row rank. It was then shown in~\cite{Willems05} that when such a PE input is applied to a controllable LTI system, the resulting input-state or input-output Hankel matrices satisfy certain rank conditions. Now known as the fundamental lemma, this result has received a great amount of attention in recent years and has been successfully used in a wide range of applications (see \cite{Markovsky21} for a comprehensive survey).\par
For certain classes of nonlinear systems, various works have proposed extensions and applications of the fundamental lemma using, e.g., basis functions. For instance, Hammerstein-Wiener systems \cite{Berberich20}, flat and feedback linearizable systems \cite{AlsaltiBerLopAll2021, Alsalti2022}, and control design of input-affine nonlinear systems \cite{DePersis22} have been considered. In these works, suitable PE conditions involving the sequence of basis functions (which depend on inputs and/or states/outputs) need to be satisfied. However, in \cite{Berberich20,AlsaltiBerLopAll2021,Alsalti2022,DePersis22} these PE conditions could only be verified a posteriori, i.e., after performing an experiment and collecting state/output data, and no a priori input design was proposed to this end.\par
In general, there exist only few results on the design of suitable inputs that result in satisfaction of the required PE conditions, specifically for nonlinear systems. In \cite{vanWaarde22}, an \emph{online} method was proposed to design inputs that result in the desired rank conditions of \cite{Willems05} on input-state or input-output Hankel matrices for linear systems. However, the resulting input is not universal in the sense that it is tailored specifically to the system on which the experiment is performed. An extension of this result to the class of bilinear systems appeared in \cite{Yuan22}. Furthermore, \cite{DePersis21} shows how suitable scaling of the initial conditions and the input to a nonlinear system leads to satisfaction of certain rank conditions on the input-state data. This enables local stabilization of an unknown nonlinear system in the first approximation but, in general, cannot be used for  the methods in \cite{Berberich20,AlsaltiBerLopAll2021,Alsalti2022,DePersis22,Alsalti2021c}.\par
The contributions of this paper are as follows. First, we propose a simple input sequence which satisfies standard PE conditions for linear systems \cite{Green86,Willems05}. Based on this result, as a second contribution, we design inputs that guarantee PE of any order for an arbitrary sequence of basis functions for the class of Hammerstein systems. The third contribution addresses nonlinear systems that are locally reachable at the origin. Specifically, we show existence of sparse input sequences that guarantee collective PE of sequences of basis functions. For single-input single-output (SISO) flat systems (which are locally reachable at the origin), we systematically design inputs that guarantee collective PE of any order of specific choices of basis functions. Finally, we illustrate the results by computing data-based controllers for SISO flat systems (as proposed in~\cite{DePersis22}).\par
Section \ref{sec_prel} introduces notation and necessary preliminaries. Sections \ref{sec_main}-\ref{sec_main3} contain the main contributions of the paper. Section \ref{sec_examples} illustrates the results with a numerical example and Section \ref{sec_conc} concludes the paper.
	\section{Notation and Preliminaries}\label{sec_prel}
\emph{Notation:} Let $\mathbb{Z}_{>0}$ denote the set of positive integers and let $\mathbb{Z}_{[a,b]}$ denote the set of integers in the interval $[a,b]$. Let the $m\times m$ identity matrix be $I_m$ and its columns be $e_i$ for $i\in\mathbb{Z}_{[1,m]}$. We use $0_{n\times m}$ to denote an $n\times m$ matrix of zeros; when the dimensions are clear from the context, we omit the subscript for notational simplicity. For a sequence $\{z_k\}_{k=0}^{N-1}$ with $z_k\in\mathbb{R}^\eta$, we denote its stacked vector as $z = \begingroup\setlength\arraycolsep{2pt}\begin{bmatrix}z_0^\top &z_1^\top & \dots & z_{N-1}^\top\end{bmatrix}\endgroup^\top$ and a stacked window of it as $z_{[l,j]} = \begingroup\setlength\arraycolsep{2pt}\begin{bmatrix}z_l^\top &z_{l+1}^\top & \dots & z_{j}^\top\end{bmatrix}\endgroup^\top$ for $0\leq l<j$. We write $M\succ0$ ($M\succeq0$) if the matrix $M$ is positive (semi-)definite.\par
In this paper, we are concerned with the notion of persistence of excitation which originated in the fields of system identification \cite{Ljung87} and adaptive control \cite{astrom08}. Several definitions of PE have appeared in the literature. The following is one of such definitions for a finite length discrete-time sequence.
\begin{definition}[\cite{Green86}]\label{def_oldPE}
		For $N\in\mathbb{Z}_{>0}$, the sequence \(\{z_k\}_{k=0}^{N-1}\), $z_k\in\mathbb{R}^{\eta}$, is exciting over the interval $[0,N-1]$ if, for some $\nu>0$, the following holds $\sum\limits_{k=0}^{N-1}z_kz_k^\top\succeq \nu I_\eta \succ0$.
\end{definition}%
In recent works on data-driven control (cf. \cite{Markovsky21}), the following definition has been commonly used.
\begin{definition}[\cite{Willems05}] For $L\in\mathbb{Z}_{>0}$ and $N\geq L$, the sequence \(\{z_k\}_{k=0}^{N-1}\), $z_k\in\mathbb{R}^{\eta}$, is persistently exciting of order \(L\) if \(\textup{rank}(\mathscr{H}_{L}(z))=\eta L\), where $\mathscr{H}_L(z) = \begin{bmatrix}
		z_{[0,L-1]} & z_{[1,L]} & \cdots & z_{[N-L,N-1]}
	\end{bmatrix}$.
	\label{def_PE}
\end{definition}

The two definitions are related as follows: Notice that if a sequence is PE of order $L$ in the sense of Definition~\ref{def_PE}, then
	\begin{equation*}
		\mathscr{H}_L(z)(\mathscr{H}_L(z))^\top = \sum\limits_{k=0}^{N-L}z_{[k,k+L-1]}z_{[k,k+L-1]}^\top\succ0.
	\end{equation*}
This means that for $L=1$, the two definitions are equivalent. If a sequence is PE in the sense of Definition~\ref{def_PE} of order $L>1$, then it is also exciting in the sense of Definition~\ref{def_oldPE}, but the converse is not necessarily true. This is because if rank$(\mathscr{H}_L(z))=\eta L$, then rank$(\mathscr{H}_{\bar{L}}(z))=\eta\bar{L}$ for any depth $\bar{L}<L$, but the converse is, in general, not true.\par
The advantage of Definition~\ref{def_PE} is that it quantifies the order of which a signal is exciting. As shown in \cite{Willems05}, this notion has an important application. In particular, if an input to a controllable LTI system is PE of order $L+n$, then the resulting input/output data matrix contains in its span any length$-L$ input/output trajectory of the system. This became known as the fundamental lemma and is summarized below.
\begin{theorem}[\cite{Willems05}]
	Let \(\{u_k,y_k\}_{k=0}^{N-1}\) be an input-output trajectory of a controllable LTI system. If \(\{u_k\}_{k=0}^{N-1}\) is PE of order $L+n$, then any \(\{\bar{u}_k,\bar{y}_k\}_{k=0}^{L-1}\) is a trajectory of the system, if and only if there exists \(\beta\in\mathbb{R}^{N-L+1}\) such that
	\begin{equation}
		\begin{bmatrix} \mathscr{H}_L(u)\\ \mathscr{H}_L(y)\end{bmatrix} \beta = \begin{bmatrix}\bar{u} \\ \bar{y}\end{bmatrix}.
		\label{eqn_fundamental_lemma}
	\end{equation}\label{thm_FL}
\vspace{-1em}
\end{theorem}
The following are recent extensions of Definition~\ref{def_PE}, which we use in our paper. The notion of \textit{collective} persistence of excitation was defined in \cite{vanWaarde20}, and extends Definition~\ref{def_PE} to multiple sequences.
\begin{definition}[\cite{vanWaarde20}]\label{def_cPE}
	For $r,L\in\mathbb{Z}_{>0}$, $j\in\mathbb{Z}_{[1,r]}$, and $N_j\geq L$, the sequences $\{z_k^{(j)}\}_{k=0}^{N_j-1}$, with $z_k^{(j)}\in\mathbb{R}^\eta$, are \textit{collectively} persistently exciting of order $L$ if rank$(\mathcal{H}_L(\mathscr{Z}))=\eta L$, where $\mathscr{Z} = \begin{bmatrix}
		(z^{(1)})^\top & \cdots & (z^{(r)})^\top
	\end{bmatrix}^\top,$ and
	\begin{equation*}
		\mathcal{H}_L(\mathscr{Z}) = \begin{bmatrix}
			\mathscr{H}_L(z^{(1)}) & \cdots & \mathscr{H}_L(z^{(r)})
		\end{bmatrix}.
	\end{equation*}
\end{definition}%

Another extension was proposed in \cite{Coulson22}, where a \textit{quantitative} notion of PE is defined in order to establish a robust version of Theorem~\ref{thm_FL}. There, a sequence $\{z_k\}_{k=0}^{N-1}$ is said to be $\alpha-$PE of order $L$ if $\sigma_{\textup{min}}(\mathscr{H}_L(z))\geq\alpha>0$, provided the Hankel matrix has at least as many columns as rows, where $\sigma_{\textup{min}}$ denotes the minimum singular value.\par
In the next section, we propose a simple input sequence which satisfies all the above definitions of PE.
	\section{PE for LTI systems}\label{sec_main}
It is well known that random inputs satisfy Definition~\ref{def_PE} with high probability. However, this is not guaranteed and in some practical applications it may not be possible/desirable to apply frequent inputs. This is the case, e.g., in networked systems with bandwidth constraints~\cite{Siami21} or closed-loop control of fluid resuscitation~\cite{XinJin19}, where sparse inputs in the form of boluses are used to identify patient-specific parameters.

In the following theorem, we propose a simple and highly sparse input sequence that is guaranteed to be PE. In particular, PE of order $L\in\mathbb{Z}_{>0}$ is achieved by giving a pulse at each of the $m$ input channels every $L$ instants.
\begin{theorem}\label{thm_PEinputLTI}
	Let $L\in\mathbb{Z}_{>0}$ and $N\geq(m+1)L-1$. If the sequence $\{u_k\}_{k=0}^{N-1}$ takes the form
	\begin{equation}
		u_k = \begin{cases}
			e_j, \qquad &k=jL-1,\\
			\mathbf{0}, \qquad &\textup{otherwise},
		\end{cases}\label{eqn_PEinputLTI}
	\end{equation}
	where $j\in\mathbb{Z}_{[1,m]}$, then it holds that $\textup{rank}(\mathscr{H}_L(u))=mL$.
\end{theorem}
\begin{proof}\let\qedsymbol\relax
	Without loss of generality, let $N=(m+1)L-1$. The following matrix has full row rank
	\begin{equation*}
		\mathscr{H}_{L}(u) = \begingroup\setlength\arraycolsep{3.5pt}\begin{bmatrix}
			\mathbf{0} & \mathbf{0} & \cdots & e_1 & \cdots & \mathbf{0} & \mathbf{0} & \cdots  & e_m\\[-1ex]
			\vdots & \vdots & \reflectbox{$\ddots$} & \vdots &\cdots & \vdots& \vdots & \reflectbox{$\ddots$} & \vdots\\[-0.5ex]
			\mathbf{0} & e_1 & \cdots & \mathbf{0} & \cdots & \mathbf{0} & e_m & \cdots & \mathbf{0}\\[-0.5ex]
			e_1 & \mathbf{0} & \cdots & \mathbf{0} & \cdots & e_m & \mathbf{0} & \cdots &\mathbf{0}
		\end{bmatrix}\endgroup.
	\end{equation*}
\end{proof}
\begin{remark}For the input in \eqref{eqn_PEinputLTI}, all singular values of $\mathscr{H}_L(u)$ are equal to one. By scaling \eqref{eqn_PEinputLTI} by $\alpha$, one obtains an input which is $\alpha-$PE of order~$L$ (cf.~\cite{Coulson22}). Note that larger values of $\alpha$ imply higher levels of PE of order~$L$.\label{remark_qPE}\end{remark}

Theorem \ref{thm_PEinputLTI} provides an explicit formula \eqref{eqn_PEinputLTI} for an input which is guaranteed to be persistently exciting in the sense of Definition \ref{def_PE} (and hence, also Definition~\ref{def_oldPE}). This can be used to apply the results of Theorem~\ref{thm_FL}. A necessary condition for the input in \eqref{eqn_PEinputLTI} to be PE of order $L+n$ is that it is at least of length $N\geq(m+1)(L+n)-1$ (cf. Definition~\ref{def_PE}). In contrast, the procedure from \cite{vanWaarde22} uses online state or output measurements to design an input such that the resulting input-output data matrix in \eqref{eqn_fundamental_lemma} has rank $mL+n$ and requires $N=(m+1)L+n-1$ samples to this end, implying that it is  sample efficient. However, the resulting input only guarantees the rank condition for the system on which the experiment was done. In contrast, the input \eqref{eqn_PEinputLTI}, although not as sample efficient, guarantees the results of Theorem~\ref{thm_FL} independently of the considered system.

Apart from its use in data-driven analysis and control, Theorem~\ref{thm_FL} also allows for identification of the LTI system's matrices up to a similarity transformation (cf. \cite[Sec. 4.4]{Willems05}). In this sense, employing \eqref{eqn_PEinputLTI} to obtain \eqref{eqn_fundamental_lemma} can be interpreted as the data-driven counterpart to classical works on identification of LTI systems from their impulse response~\cite{HO66}.

In Section \ref{sec_main2}, we consider Hammerstein nonlinear systems and show how to guarantee PE of basis functions that depend on the input, by the design of the input only.
	\section{PE for Hammerstein systems}\label{sec_main2}
Consider an unknown Hammerstein system of the form
\begin{equation}
	\begin{aligned}
		x_{k+1} = Ax_k + B\gamma(u_k), \quad y_k = Cx_k + D\gamma(u_k),
	\end{aligned}\label{eqn_HamSys}
\end{equation}
with $(A,B)$ controllable, $x_k\in\mathbb{R}^n,\,u_k\in\mathbb{R}^m,$ and $y_k\in\mathbb{R}^p$. Furthermore, $\gamma:\mathbb{R}^m\to\mathbb{R}^{\bar{m}}$ is an unknown nonlinear function which satisfies $\gamma(\mathbf{0})=\mathbf{0}$ without loss of generality. An extension of Theorem \ref{thm_FL} to the class of Hammerstein systems appeared in \cite{Berberich20}, assuming that $\gamma_i$, $i\in\mathbb{Z}_{[1,\bar{m}]}$, belong to the span of a given set of $r$ basis functions $\psi_j:\mathbb{R}^m\to\mathbb{R}$, $j\in\mathbb{Z}_{[1,r]}$, which satisfy $\psi_j(\mathbf{0})=0$ (note that functions with $\psi_j(\mathbf{0})\neq0$ can be suitably shifted by a constant). In this case, PE must be imposed on the sequence of basis functions. In this section, we show how this can be done only by the design of the physical input $u$.\par
We denote the stacked vector of the basis functions by $\Psi(u_k)=\begingroup
\setlength\arraycolsep{1pt}\begin{bmatrix} 	\psi_1(u_k) & \cdots & \psi_r(u_k) \end{bmatrix}^\top\endgroup$. For $\lambda_j\in\mathbb{R}^m$, $j\in\mathbb{Z}_{[1,r]}$, we define the following square matrix $\Lambda\in\mathbb{R}^{r\times r}$
\begin{equation}
	\begin{aligned}
		\Lambda &= \begin{bmatrix}
			\Psi(\lambda_1) & \Psi(\lambda_2) & \cdots & \Psi(\lambda_r)
		\end{bmatrix}.
	\end{aligned}\label{eqn_Lambda}
\end{equation}
The following theorem provides conditions for an input sequence $\{u_k\}_{k=0}^{N-1}$ that is guaranteed to result in a PE sequence of basis functions $\{\hat{\Psi}_k\}_{k=0}^{N-1}$, where $\hat{\Psi}_k\coloneqq\Psi(u_k)$.
\begin{theorem}\label{thm_PEinputHam}
	Given $L\in\mathbb{Z}_{>0}$ and any $r$ linearly independent basis functions $\psi_j:\mathbb{R}^m\to\mathbb{R}$ which satisfy $\psi_j(\mathbf{0})=0$, let $N\geq(r+1)L-1$. Let $\{u_k\}_{k=0}^{N-1}$ take the form
	\begin{equation}
		u_k = \begin{cases}
			\lambda_j, \qquad &k=jL-1,\\
			\mathbf{0}, \qquad &\textup{otherwise},
		\end{cases}\label{eqn_PEinputHam}
	\end{equation}
	where $j\in\mathbb{Z}_{[1,r]}$. If $\lambda_j\in\mathbb{R}^m$ are such that the matrix $\Lambda$ in \eqref{eqn_Lambda} is invertible, then it holds that $\textup{rank}(\mathscr{H}_L(\hat{\Psi})) = rL.$
\end{theorem}
\begin{proof}
	Let $\{\tilde{u}_k\}_{k=0}^{N-1}$, $\tilde{u}_k\in\mathbb{R}^r$, take the form of \eqref{eqn_PEinputLTI}. Then by Theorem~\ref{thm_PEinputLTI} it holds that rank$(\mathscr{H}_L(\tilde{u}))=rL$. Next, notice that $\mathscr{H}_L(\hat{\Psi}) = (I_L\otimes \Lambda)\mathscr{H}_L(\tilde{u})$, where $\otimes$ denotes the Kronecker product. Since $\Lambda$ is invertible by assumption, it follows that $(I_L\otimes\Lambda)$ is a square full rank matrix and hence rank$(\mathscr{H}_L(\hat{\Psi}))=\textup{rank}(\mathscr{H}_L(\tilde{u}))=rL$.
\end{proof}

In order to find the values of $\lambda_j$ which make $\Lambda$ invertible, one can formulate a feasibility problem as follows
\begin{equation}
	\begin{aligned}
		\textup{find}\quad\lambda_1,\ldots,\lambda_r, \qquad
		\textup{s.t.}\quad\textup{rank}(\Lambda)=r.
	\end{aligned}\label{eqn_feasibility}
\end{equation}
This is equivalent to requiring that $\Lambda^\top \Lambda\succ0$. Such a problem can be solved offline, e.g., iteratively as in \cite{Sun17} or by reformulating it as a regularized unconstrained nonlinear least squares problem \cite{Markovsky13_SLRA}. Solving \eqref{eqn_feasibility} only requires knowledge of the user-defined basis functions. The following theorem shows that for \textit{any} given set of linearly independent basis functions, \eqref{eqn_feasibility} is feasible.
\begin{theorem}\label{thm_alwaysexistslambda}
	Given any set of $r$ linearly independent basis functions $\psi_j:\mathbb{R}^m\to\mathbb{R}$, there exist $\lambda_j\in\mathbb{R}^m$, for $j\in\mathbb{Z}_{[1,r]}$, such that $\Lambda$ in \eqref{eqn_Lambda} is invertible.
\end{theorem}
\begin{proof}
	Suppose by contradiction that the maximum rank that $\Lambda$ can attain for arbitrary $\lambda_j \in \mathbb{R}^m$, $j\in\mathbb{Z}_{[1,r]}$, is $d<r$. Consider such a choice of $\lambda_j \in \mathbb{R}^m$, $j\in\mathbb{Z}_{[1,r]}$, which results in rank$(\Lambda)=d<r$. Then, there exists a non-zero vector $\rho~\in~\mathbb{R}^r$ such that $\rho^\top \Lambda = 0$. Moreover, for any $\bar\lambda \in \mathbb{R}^m$, the matrix $\begin{bmatrix}\Lambda & \Psi(\bar\lambda)\end{bmatrix}$ has the same rank as $\Lambda$. Hence, it holds that $\rho^\top\Psi(\bar\lambda)=0$ for arbitrary $\bar\lambda \in \mathbb{R}^m$, i.e.,
		\begin{equation}
			\rho_1\psi_1(\bar\lambda) + \cdots + \rho_r\psi_r(\bar\lambda) =0,\qquad \forall \bar\lambda \in \mathbb{R}^m.\label{eqn_LIbasisfcns}
		\end{equation}
		This, however, contradicts linear independence of the basis functions, thus proving that there exists $\lambda_j\in\mathbb{R}^m$, $j\in\mathbb{Z}_{[1,r]}$, such that $\Lambda$ is invertible.
\end{proof} 

In the next section, we study the class of locally reachable nonlinear systems at the origin and show existence of sparse inputs that guarantee collective PE of basis functions. For SISO flat systems (which are locally reachable), we show how to design those inputs such that PE of any order $L>0$ can be guaranteed for a specific choice of basis functions.
	\section{PE for locally reachable nonlinear systems}\label{sec_main3}
Consider an unknown nonlinear system of the form
\begin{equation}
	x_{k+1} = f(x_k,u_k),\label{eqn_NLsys}
\end{equation}%
with $x_k\in\mathbb{R}^n,u_k\in\mathbb{R}^m$ being the state and input vectors, respectively, and $f:\mathbb{R}^n\times\mathbb{R}^m\to\mathbb{R}^n$ is an unknown function satisfying $f(\mathbf{0},\mathbf{0})=\mathbf{0}$. For $\mu\in\mathbb{Z}_{>0}$, the set of all states which can be reached from $x_0$ in $\mu$ steps is defined as
\begin{align}
	&\mathcal{R}_\mu(x_0) \hspace{-1mm}=\hspace{-1mm} \left\lbrace \hspace{-0.5mm} x_\mu\in\mathbb{R}^n \,\left|\, \begin{aligned}
		&\exists\, u_{[0,\mu-1]},\, u_k\in\mathbb{R}^m,\\ &\textup{s.t.}\,\,  x_{k+1}\hspace{-0.5mm}=\hspace{-0.5mm}f(x_k,u_k),\,\forall k\hspace{-0.5mm}\in\hspace{-0.5mm}\mathbb{Z}_{[0,\mu-1]}.
	\end{aligned} \hspace{-0.5mm} \right. \right\rbrace \hspace{-1mm}.\notag
\end{align}
It was shown in \cite{Melkior20} how one can, under certain assumptions, obtain a guaranteed under-approximation of the reachable set of a nonlinear system with unknown dynamics. For the remainder of this section, we make the following assumption.
\begin{assumption}\label{assmp_reachability}
	For the system \eqref{eqn_NLsys}, there exists $\mu\in\mathbb{Z}_{>0}$ such that the origin is contained in the interior of $\mathcal{R}_\mu(\mathbf{0})$.
\end{assumption}%

Assumption~\ref{assmp_reachability} implies that the system is locally reachable at the origin. A sufficient condition for local reachability at $x=\mathbf{0}$ is that the linearization of system~\eqref{eqn_NLsys} at the origin is controllable (cf. \cite[Lemma 3.7.8]{Sontag13}).

Consider now $r$ basis functions $\theta_j:~\mathbb{R}^n\times \mathbb{R}^m\to\mathbb{R}$ which satisfy $\theta_j(\mathbf{0},\mathbf{0})=0$, $j\in\mathbb{Z}_{[1,r]}$, and denote the stacked vector of them by $\Theta(x_k,u_k)=\begingroup
\setlength\arraycolsep{1pt}\begin{bmatrix}\theta_1(x_k,u_k) & \cdots & \theta_r(x_k,u_k) \end{bmatrix}^\top\endgroup$. Suppose that the functions are linearly independent on arbitrary domains with non-empty interior\footnote{This assumption is satisfied for sinusoidal functions, exponential functions and monomials, among others \cite{Christensen06}. Note that basis functions that do not satisfy $\theta_j(\mathbf{0},\mathbf{0})=0$ can be suitably shifted by a constant.} $D_x\times D_u\subset\mathbb{R}^n\times~\mathbb{R}^m$. The objective is to design inputs $\{u_k^{(j)}\}_{k=0}^{N_j-1}$ such that the sequences of basis functions $\{\hat{\Theta}_k^{(j)}\}_{k=0}^{N_j-1}$ (with $\hat{\Theta}_k^{(j)}\coloneqq\Theta(x_k^{(j)},u_k^{(j)})$) are collectively persistently exciting of order $L$. According to Definition \ref{def_cPE}, the following mosaic Hankel matrix must have full row rank
\begin{align}
	&\mathcal{H}_L(\vartheta) =\label{eqn_MosaicHankelTheta} \begin{bmatrix}
		\mathscr{H}_L(\hat{\Theta}^{(1)}) & \cdots & \mathscr{H}_L(\hat{\Theta}^{(r)})
	\end{bmatrix},
\end{align}
where $\vartheta = \begin{bmatrix}
	(\hat{\Theta}^{(1)})^\top & \cdots & (\hat{\Theta}^{(r)})^\top
\end{bmatrix}^\top$.

Consider a submatrix of \eqref{eqn_MosaicHankelTheta} composed of the $L+\mu$ element of each of the $r$ sequences of basis functions
\begin{align}
	&W =
	\begin{bmatrix}
		\hat{\Theta}^{(1)}_{L+\mu-1} & \hat{\Theta}^{(2)}_{L+\mu-1} & \cdots & \hat{\Theta}^{(r)}_{L+\mu-1}
	\end{bmatrix}.\label{eqn_W}
\end{align}
Similar to Section \ref{sec_main2}, we would like $W$ to be invertible. Since the state at time $x_{L+\mu-1}^{(j)}$ is not a free variable, ensuring invertibility of $W$ is a control problem. In particular, one must select the inputs such that the corresponding state and input pairs at time $L+\mu-1$ result in an invertible matrix $W$. In a data-driven setting, such a control problem is difficult to solve since no model knowledge is available and - to begin with - no persistently exciting data is available yet to apply data-driven control techniques. Therefore, we first show in Lemma~\ref{lemma_rexperimentsexist} that there exist such $r$ input sequences $u_{[0,L+\mu-1]}^{(j)}$, $j\in\mathbb{Z}_{[1,r]}$, which make $W$ invertible and then prove in Theorem~\ref{thm_generalPE} how invertibility of $W$ results in collective PE of the basis functions. Later, in the following Subsection~\ref{sec_flatPE} we illustrate how, under suitable assumptions, one can indeed find the desired control inputs a priori, which guarantee collective PE of any order for sequences of basis functions that depend on input and output data of SISO flat systems.

\begin{lemma}\label{lemma_rexperimentsexist}
	Let Assumption \ref{assmp_reachability} be satisfied and suppose that $r$ basis functions $\theta_j:\mathbb{R}^n\times\mathbb{R}^m\to\mathbb{R}$ are linearly independent on $\mathcal{R}_\mu(\mathbf{0})\times\mathbb{R}^m$. Then there exist $r$ sequences $u_{[0,L+\mu-1]}^{(j)}$, $j\in\mathbb{Z}_{[1,r]}$, which when applied to \eqref{eqn_NLsys} starting from~$x_0^{(j)}=\mathbf{0}$, result in $x_{L+\mu-1}^{(j)}$, such that $W$ in \eqref{eqn_W} is invertible.
\end{lemma}
\begin{proof}
	Using similar arguments to the proof of Theorem~\ref{thm_alwaysexistslambda}, it can be shown by linear independence of the basis functions $\theta_j$ on $\mathcal{R}_\mu(\mathbf{0})\times\mathbb{R}^m$ that there exists $r$ pairs $(x_{L+\mu-1}^{(j)},u_{L+\mu-1}^{(j)})\in \mathcal{R}_\mu(\mathbf{0})\times\mathbb{R}^m$ such that $W$ is invertible. 
	
	Since $f(\mathbf{0},\mathbf{0})=\mathbf{0}$, then starting from zero initial conditions and setting the input to $u^{(j)}_{[0,L-2]}=0$, one can express the state $x_{L+\mu-1}^{(j)}$ in terms of $u^{(j)}_{[L-1,L+\mu-2]}$ only, i.e., $x_{L+\mu-1}^{(j)} = f(f\cdots(f(\mathbf{0},u_{L-1}^{(j)}),u_{L}^{(j)})\cdots,u_{L+\mu-2}^{(j)})$. Finally, since the system is locally reachable by Assumption~\ref{assmp_reachability}, there exist inputs $u_{[L-1,L+\mu-1]}^{(j)}$, $j\in\mathbb{Z}_{[1,r]}$, which steer the system from $x_{L-1}^{(j)}=\mathbf{0}$ to $x_{L+\mu-1}^{(j)}$ in $\mu$~steps.
\end{proof}

The following theorem shows how to use the results of Lemma \ref{lemma_rexperimentsexist} to obtain collectively persistently exciting sequences of basis functions of any order $L$.
\begin{theorem}\label{thm_generalPE}
	Let Assumption~\ref{assmp_reachability} hold. Given $L,\mu\in\mathbb{Z}_{>0}$ and $r$ basis functions $\theta_j:\mathbb{R}^n\times\mathbb{R}^m\to\mathbb{R}$ that are linearly independent on $\mathcal{R}_\mu(\mathbf{0})\times\mathbb{R}^m$ and satisfy $\theta_j(\mathbf{0},\mathbf{0})=0$, let $N_j\geq2L+\mu-1$ for $j\in\mathbb{Z}_{[1,r]}$. Furthermore, let the sequences $\{u_k^{(j)}\}_{k=0}^{N_j-1}$ take the form
	\begin{equation}
		u_k^{(j)} = \begin{cases}
			\eta^{(j)}_{[0,\mu]}, \quad & k\in\mathbb{Z}_{[L-1,L+\mu-1]},\\
			\mathbf{0}, \quad &\textup{otherwise},
		\end{cases}\label{eqn_PElocallyreachable}
	\end{equation}
	where $j\in\mathbb{Z}_{[1,r]}$, and $\eta^{(j)}_{[0,\mu]}$ are such that $W$ in \eqref{eqn_W} is invertible. If \eqref{eqn_PElocallyreachable} are applied to \eqref{eqn_NLsys} starting from $x_0^{(j)}=\mathbf{0}$, then for the matrix in \eqref{eqn_MosaicHankelTheta} it holds that $\textup{rank}(\mathcal{H}_L(\vartheta)) = rL$.
\end{theorem}
\begin{proof}
	Each block row of \eqref{eqn_MosaicHankelTheta} has $r$ linearly independent columns given by the columns of $W$. Notice that each $\mathscr{H}_L(\hat{\Theta}^{(j)})$, $j\in\mathbb{Z}_{[1,r]}$, in \eqref{eqn_MosaicHankelTheta} has a lower block-anti-triangular structure due to $f(\mathbf{0},\mathbf{0})=\mathbf{0},\,\Theta(\mathbf{0},\mathbf{0})=\mathbf{0}$ and the choice of the inputs \eqref{eqn_PElocallyreachable}. As a result, every block row \eqref{eqn_MosaicHankelTheta} is linearly independent from the others. Since there are $L$ such block rows, it holds that rank$(\mathcal{H}_L(\vartheta))=rL$.
\end{proof}

Notice that the results of Lemma~\ref{lemma_rexperimentsexist} for locally reachable nonlinear systems are analogous to that of Theorem~\ref{thm_alwaysexistslambda} for Hammerstein systems. However, formulating a nonlinear feasibility problem similar to \eqref{eqn_feasibility} to find $u_{[0,L+\mu-1]}^{(j)}$, $j\in\mathbb{Z}_{[1,r]}$ would require knowledge of the unknown function $f$.

It was observed in simulations that randomly sampling the input sequences $\eta^{(j)}_{[0,\mu]}$, $j\in\mathbb{Z}_{[1,r]}$, in Theorem~\ref{thm_generalPE} from a uniform distribution typically results in a corresponding invertible matrix $W$. However, such a heuristic approach is not always guaranteed to achieve this result. To systematically find the desired input sequences, one must impose additional assumptions on the {class} of systems and the {choice} of basis functions. To this end, we consider in the next subsection SISO flat nonlinear systems (which are locally reachable at the origin), and show how one can guarantee PE of any order $L>0$ a priori, for a specific choice of basis functions.
\subsection{SISO flat nonlinear systems}\label{sec_flatPE}
Consider an unknown SISO flat system of the form
\begin{equation}
	\begin{aligned}
		x_{k+1} = f(x_k,u_k), \quad y_k = h(x_k),
	\end{aligned}\label{eqn_flatsys}
\end{equation}
where \(x_k\in\mathbb{R}^n, u_k,\,y_k\in\mathbb{R}\) and $f:\mathbb{R}^n\times\mathbb{R}\to\mathbb{R}^n$, $h:~\mathbb{R}^n\to~\mathbb{R}$ are smooth unknown functions with $f(\mathbf{0},0)=\mathbf{0}$ and $h(\mathbf{0})=0$. Let $f_O^j(x_k)$ denote the $j-$th iterated composition of the undriven dynamics $f(x_k,0)$.\par
Since the system is flat (i.e., has a well defined relative degree equal to the system dimension $n$, cf. \cite[Sec. III.A]{AlsaltiBerLopAll2021}), it can be transformed into the discrete-time normal form provided that $0\in\textup{Im}\left(h(f_O^{n-1}(f(x,\cdot)))\right)$ holds for all $x\in\mathbb{R}^n$ (cf. \cite[Sec. 2]{MonacoNor1987} for more details). This means that there exists an invertible (w.r.t. ${v}_k$) control law ${u}_k=q({x}_k,{v}_k)$, with $q:\mathbb{R}^n\times\mathbb{R}\to\mathbb{R}$ and an invertible coordinate transformation $\xi_k=T({x}_k)=y_{[k,k+n-1]}$, such that
\begin{equation}
	\begin{matrix}
		\xi_{k+1} = {A}\xi_k + {B}{v}_k, \qquad
		{y}_k = {C}\xi_k,
	\end{matrix}
	\label{BINF}%
\end{equation}%
Furthermore, $A,B,C$ are in the Brunovsky canonical form (cf. \cite[Thm.~2]{AlsaltiBerLopAll2021}) which is a controllable/observable triplet. Hence, the system is $n$ steps locally reachable at the origin. The synthetic input $v_k$ takes the form
\begin{equation}
	\begin{aligned}
		v_k = h(f_O^{n-1}(f(x_k,u_k))).
	\end{aligned}\label{eqn_expressionforv}
\end{equation}

A sufficient condition for the analogue of Theorem~\ref{thm_FL} to flat systems \cite[Prop. 1]{AlsaltiBerLopAll2021}, and for designing controllers from data in \cite[Cor. 2]{DePersis22}, is persistence of excitation of a sequence of basis functions which contain $h\circ f_O^{n-1}\circ f$ in their span. To check the PE condition, one typically performs an experiment of length $N\geq (r+1)L-1$, collects the corresponding state or output measurements and then verifies the rank of the resulting Hankel matrix.

In this section, we illustrate how one can enforce PE of any order for a \textit{particular choice} of basis functions \textit{a priori}. A specific choice of basis functions may, in general, not contain the unknown nonlinearity \eqref{eqn_expressionforv} in its span. Nonetheless, enforcing PE of such basis functions is still useful for, e.g., designing locally stabilizing controllers for unknown SISO flat systems (cf. \cite[Sec. VII.B]{DePersis22}), and for data-driven nonlinear predictive control \cite{Alsalti2021c}, provided that the basis functions result in a good local approximation of~\eqref{eqn_expressionforv}. 

Since the map from $u$ to $v$ is invertible and since $f(\mathbf{0},0)=\mathbf{0}$ and $h(\mathbf{0})=0$, a non-zero input applied to the system from zero initial conditions results in a non-zero value of $v$ in \eqref{eqn_expressionforv}. Moreover, invertibility implies that for all $\delta_1,\delta_2\in\mathbb{R}$, the following holds
\begin{equation}
	\delta_1\hspace{-0.5mm}\neq \hspace{-0.5mm}\delta_2 \hspace{-1mm}\iff\hspace{-1mm} h(f_O^{n-1}(f(\mathbf{0},\delta_{1})))\hspace{-0.5mm}\neq\hspace{-0.5mm} h(f_O^{n-1}(f(\mathbf{0},\delta_{2}))).\label{eqn_uniquenessofv}
\end{equation}

We exploit this fact to prove the following lemma, which will be needed later for the main result of this subsection.
\begin{lemma}\label{lemma_vandermonde}
	For $t\in\mathbb{Z}_{>0}$ let $\delta_j\neq0$, $j\in\mathbb{Z}_{[1,t]}$, be mutually distinct values and define $v_{\delta_j}\coloneqq h(f_O^{n-1}(f(\mathbf{0},\delta_j)))$. Then, the following matrix is invertible:
\end{lemma}
\begin{equation}
	\Omega = \begin{bmatrix}
		v_{\delta_1} & v_{\delta_2} & \cdots & v_{\delta_t}\\
		v_{\delta_1}^2 & v_{\delta_2}^2 & \cdots & v_{\delta_t}^2\\
		\vdots & \vdots & \ddots & \vdots\\
		v_{\delta_1}^{t} & v_{\delta_2}^{t} & \cdots & v_{\delta_t}^{t}
	\end{bmatrix}.\label{eqn_OmegaVandermonde}
\end{equation}
\begin{proof}
	Since for $j\in\mathbb{Z}_{[1,t]}$, $\delta_j\neq0$ are mutually distinct values, it follows that the corresponding values $v_{\delta_j}$ are also distinct and non-zero (compare the discussion above Lemma~\ref{lemma_vandermonde}). The matrix $\Omega$ can be written as $\Omega = V^\top \Delta$, where $V\in\mathbb{R}^{t\times t}$ is a square Vandermonde matrix composed of the distinct $v_{\delta_j}$ and, hence, invertible and $\Delta\in\mathbb{R}^{t\times t}$ is a diagonal matrix containing $v_{\delta_j}$. The proof is concluded by noting that $V$ and $\Delta$ are invertible matrices.
\end{proof}

In the following, we consider monomial basis functions in the transformed state and input up to some finite order $t\in\mathbb{Z}_{>0}$, and hence $r=t(n+1)$.
\begin{align}
	&\Theta(\xi_k,u_k) = \label{eqn_specificchoice}\begingroup
	\setlength\arraycolsep{2pt}\begin{bmatrix}
		u_k & u_k^2 & \cdots & u_k^t& \xi_k^\top & (\xi_k^2)^\top & \cdots & (\xi_k^t)^\top 
	\end{bmatrix}^\top\endgroup\hspace{-2mm}.
\end{align}
The powers are defined element-wise, i.e., $\xi_k^t=[\xi_{1,k}^t \,\, \cdots \,\, \xi_{n,k}^t]^\top$. Notice that these functions depend only on the inputs and outputs of \eqref{eqn_flatsys} since $\xi_k=y_{[k,k+n-1]}$ (cf. \eqref{BINF}). In the following theorem, we show how to choose input sequences $\{u_k^{(j)}\}_{k=0}^{N_j-1}$, $j\in\mathbb{Z}_{[1,r]}$, such that the resulting sequences of basis functions $\{\hat{\Theta}_k^{(j)}\}_{k=0}^{N_j-1}$ are collectively persistently exciting of order $L>0$, i.e., that the corresponding mosaic Hankel matrix $\mathcal{H}_L(\vartheta)$ of the form \eqref{eqn_MosaicHankelTheta} has full row rank.
\begin{theorem}\label{thm_aprioriFL}
		For $t\in\mathbb{Z}_{>0}$, let $\delta_j\neq0$, $j\in~\mathbb{Z}_{[1,t(n+1)]}$, be mutually distinct values. For $L\in\mathbb{Z}_{>0}$, $N_j\geq 2L+n-1$ and the basis functions in \eqref{eqn_specificchoice}, let $\{u_k^{(j)}\}_{k=0}^{N_j-1}$ take the form in \eqref{eqn_PElocallyreachable} with the corresponding $\eta^{(j)}_{[0,n]}$ given by
		\begin{equation}
			\eta^{(j)}_{[0,n]} = \begin{cases}
				\begin{bsmallmatrix}
					\mathbf{0}_{n-j\times1}\\ \delta_j\\ \mathbf{0}_{j\times1}
				\end{bsmallmatrix}, \quad & \textup{for }j\in\mathbb{Z}_{[1,n]},\\
				&\boldsymbol{\vdots}\\
				\begin{bsmallmatrix}
					\mathbf{0}_{tn-j\times1}\\ \delta_j\\ \mathbf{0}_{j-(t-1)n\times1}
				\end{bsmallmatrix}, \quad & \textup{for }j\in\mathbb{Z}_{[(t-1)n,tn]},\\
				\begin{bsmallmatrix}
					\mathbf{0}_{n\times1}\\ \delta_j
				\end{bsmallmatrix}, \quad & \textup{for }j\in\mathbb{Z}_{[tn+1,t(n+1)]}.
			\end{cases}
			\label{eqn_PEinputSISOflat}
		\end{equation}
		If \eqref{eqn_PEinputSISOflat} are applied to \eqref{eqn_flatsys} starting from zero initial conditions, then $\textup{rank}(\mathcal{H}_L(\vartheta))=t(n+1)L$.
\end{theorem}
\begin{proof}
	Without loss of generality, let $N_j=2L+n-1$ for all $j\in\mathbb{Z}_{[1,t(n+1)]}$. For $c\in\mathbb{Z}_{[0,n-1]}$, we define $v_{\delta_j,[0,c]}~\coloneqq~\begingroup
		\setlength\arraycolsep{1pt}\begin{bmatrix}
			h(f_O^{n-1}(f(\mathbf{0},\delta_{j}))) \,\, \cdots \,\, h(f_O^{n+c-1}(f(\mathbf{0},\delta_{j})))
		\end{bmatrix}^\top\endgroup\hspace{-1.5mm}$ (with some abuse of notation we also use $v_{\delta_j,0}=v_{\delta_j}$). Since the system \eqref{BINF} is in the Brunovsky form, applying the inputs $\{u_k^{(j)}\}_{k=0}^{N_j-1}$ as defined in the theorem statement from zero initial conditions results in
		\begin{equation}
			\xi_{L+n-1}^{(j)}\hspace{-1mm} = \hspace{-1mm}\begin{cases}
				\begin{bsmallmatrix}
					\mathbf{0}_{n-j\times 1}\\ v_{\delta_j,[0,j-1]}
				\end{bsmallmatrix}, \, & \textup{for }j\in\mathbb{Z}_{[1,n]},\\
				&\boldsymbol{\vdots}\\
				\begin{bsmallmatrix}
					\mathbf{0}_{tn-j\times 1}\\ v_{\delta_j,[0,j-(t-1)n-1]}
				\end{bsmallmatrix}, \, & \textup{for }j\in\mathbb{Z}_{[(t-1)n+1,tn]},\\
				\mathbf{0},\,&\textup{for }j\in\mathbb{Z}_{[tn+1,t(n+1)]}.
			\end{cases}\label{eqn_XiStateForW}
	\end{equation}%
	Now, we consider a submatrix of $\mathcal{H}_L(\vartheta)$ of the form of $W$ in \eqref{eqn_W}. For the choice of basis functions in \eqref{eqn_specificchoice}, the inputs $\{u_k^{(j)}\}_{k=0}^{N_j-1}$ as defined in the theorem statement and the corresponding state values in \eqref{eqn_XiStateForW}, the matrix $W$ takes the form \eqref{eqn_Wflat} (see next page). Following similar arguments to the proof of Lemma~\ref{lemma_vandermonde}, one can show that $W_u\in\mathbb{R}^{t\times t}$ is invertible since $\delta_j$, $j\in\mathbb{Z}_{[tn+1,t(n+1)]}$, are non-zero and mutually distinct values.
		\begin{figure*}[!t]
			\normalsize
			\begin{align}
					W = \begin{bmatrix}
						\mathbf{0}& W_u\\
						W_\xi & \mathbf{0}
					\end{bmatrix} = \left[\begin{array}{cccc|ccc}
						0 & 0 & \cdots & 0 & \delta_{tn+1} & \cdots & \delta_{t(n+1)}\\[-1ex]
						\vdots & \vdots & \vdots & \vdots & \vdots & \vdots & \vdots\\[-0.5ex]
						0 & 0 & \cdots & 0 & \delta_{tn+1}^t & \cdots & \delta_{t(n+1)}^t\\[-0.25ex]
						\hline
						\begin{pmatrix}
							\mathbf{0}_{n-1\times1}\\[-0.5ex] v_{\delta_1}
						\end{pmatrix} & \begin{pmatrix}
							\mathbf{0}_{n-2\times1}\\[-0.5ex] v_{\delta_2,[0,1]}
						\end{pmatrix} & \cdots & \begin{pmatrix}
							v_{\delta_{tn},[0,n-1]}
						\end{pmatrix} & \mathbf{0} & \cdots & \mathbf{0}\\[-1ex]
						\vdots & \vdots & \vdots & \vdots & \vdots & \vdots & \vdots\\[-1ex]
						\begin{pmatrix}
							\mathbf{0}_{n-1\times1}\\[-0.5ex] v_{\delta_1}
						\end{pmatrix}^t & \begin{pmatrix}
							\mathbf{0}_{n-2\times1}\\[-0.5ex] v_{\delta_2,[0,1]}
						\end{pmatrix}^t & \cdots & \begin{pmatrix}
							v_{\delta_{tn},[0,n-1]}
						\end{pmatrix}^t & \mathbf{0} & \cdots & \mathbf{0}
					\end{array}
					\right].
					\label{eqn_Wflat}
				\end{align}
			\hrulefill
			\vspace{-1em}
		\end{figure*}
		Using the columns of $W_\xi\in\mathbb{R}^{tn\times tn}$ in \eqref{eqn_Wflat}, we construct $n$ submatrices $\overline{W}_{i,\xi}\in\mathbb{R}^{tn\times t}$, $i\in\mathbb{Z}_{[1,n]}$, of the form
		\begin{align*}
			&\overline{W}_{i,\xi} =\\
			&\begingroup
			\setlength\arraycolsep{1.5pt}\begin{bmatrix}
				\begin{pmatrix}
					\mathbf{0}_{n-i\times1}\\ v_{\delta_i,[0,i-1]}
				\end{pmatrix} & \begin{pmatrix}
					\mathbf{0}_{n-i\times1}\\ v_{\delta_{i+n},[0,i-1]}
				\end{pmatrix} & \cdots & \begin{pmatrix}
					\mathbf{0}_{n-i\times1}\\ v_{\delta_{i+(t-1)n},[0,i-1]}
				\end{pmatrix}\\
				\vdots & \vdots & \vdots & \vdots\\
				\begin{pmatrix}
					\mathbf{0}_{n-i\times1}\\ v_{\delta_{i},[0,i-1]}
				\end{pmatrix}^{\hspace{-0.25mm}t} & \begin{pmatrix}
					\mathbf{0}_{n-i\times1}\\ v_{\delta_{i+n},[0,i-1]}
				\end{pmatrix}^{\hspace{-0.25mm}t} & \cdots & \begin{pmatrix}
					\mathbf{0}_{n-i\times1}\\ v_{\delta_{i+(t-1)n},[0,i-1]}
				\end{pmatrix}^{\hspace{-0.25mm}t}
			\end{bmatrix}\endgroup\hspace{-1mm}.
		\end{align*}%
		Each matrix of this form has $t$ rows of the form $\Omega$ in \eqref{eqn_OmegaVandermonde}. Since $\delta_j$, $j\in\mathbb{Z}_{[1,tn]}$, are non-zero and distinct, it follows from Lemma \ref{lemma_vandermonde} that the corresponding $\Omega$ is invertible and hence, each matrix $\overline{W}_{i,\xi}$ has rank~$t$. Notice that the columns of each matrix $\overline{W}_{i,\xi}$ are linearly independent with respect to the columns of any other  $\overline{W}_{j,\xi}$, $i\neq j\in\mathbb{Z}_{[1,n]}$. This follows since (i) the rows of the form $\Omega$ appear in different rows in each $\overline{W}_{i,\xi}$ and (ii) the structure in which the block rows $\mathbf{0}_{n-i\times1}$ appear in each submatrix. As a result, rank$(W_\xi)=tn$. Due to the structure of $W_u$ and $W_\xi$, it follows that rank$(W)\hspace{-0.25mm}=\hspace{-0.25mm}t(n\hspace{-0.1mm}+\hspace{-0.1mm}1)$ and hence $W$ is invertible. Finally, it follows from Theorem~\ref{thm_generalPE} that rank$(\mathcal{H}_L(\vartheta))\hspace{-0.25mm}=\hspace{-0.25mm}t(n\hspace{-0.1mm}+\hspace{-0.1mm}1)L$.
\end{proof}

In the following section, we illustrate the results of Theorem~\ref{thm_aprioriFL} with an example.
	\section{Numerical Example}\label{sec_examples}
Consider a second order SISO flat system of the form in \eqref{BINF}, with $v_k = -\sin(x_{1,k}) + x_{1,k}x_{2,k}^2 - x_{1,k}^3x_{2,k} + u_k.$ In this example, we compare the performance of three nonlinear controllers: (i) An exact linearizing and stabilizing controller designed using basis functions that include $v$ in their span \cite[Cor. 2]{DePersis22}, and two locally stabilizing controllers (ii and iii) designed using the following choice of basis functions\footnote{The method described in \cite[Cor. 2]{DePersis22} requires that the unknown map \eqref{eqn_expressionforv} is linear in $u$, which is why we use the basis functions \eqref{eqn_ex_basisfunctions}. Although the choice of the basis functions in \eqref{eqn_ex_basisfunctions} is different from that in \eqref{eqn_specificchoice}, one can easily see from the proof of Theorem \ref{thm_aprioriFL} that using inputs of the form \eqref{eqn_PEinputSISOflat} also guarantees collective PE of \eqref{eqn_ex_basisfunctions}.} which do not contain $v$ in their span \cite[Cor. 2 and Sec.~III.B]{DePersis22}
\begin{equation}
	\Theta(\xi_k,u_k) = \begin{bmatrix}
		u_k & \xi_k^\top & (\xi_k^2)^\top & (\xi_k^3)^\top
	\end{bmatrix}^{\hspace{-0.5mm}\top}.\label{eqn_ex_basisfunctions}
\end{equation}
For all three controllers, PE of the basis functions of order one is a necessary and sufficient condition for the feasibility of the convex program that is solved to obtain the control gains (cf. \cite[Cor. 2, Thm. 2, and Thm. 5]{DePersis22}). For controllers (i) and (ii), PE is enforced by sampling the input randomly. For controller (iii), PE is enforced \textit{a priori} using the results of Theorem~\ref{thm_aprioriFL}. In this case, we used a straightforward extension of \cite[Cor. 2]{DePersis22} such that collected data from multiple experiments (i.e., collective PE) can be used to design the controller.

Since the system is unstable, the input data (of length $N=21$) for controllers (i) and (ii) had to be sampled from the uniform distribution $U(-0.25,0.25)$, whereas using multiple experiments as in Theorem~\ref{thm_aprioriFL} allowed us to use inputs (each of length $N_j=3$) with larger magnitudes (sampled from $U(-1,1)$). In \cite{vanWaarde20}, a similar observation was made for linear systems. As a result, a larger quantitative level of PE was attained (cf. Remark \ref{remark_qPE} and Table~\ref{table_comparison}).

The performance of the closed-loop system (over $T=20$ time instants) was compared starting from the same initial conditions (randomly sampled from $U(-1,1)\times U(-1,1)$). Table~\ref{table_comparison} shows the average cumulative stabilization errors (defined as $\sum_{k=0}^{T-1}\frac{1}{T}|x_{i,k}|$, for $i=1,2,\,T=20$) for all three controllers over 100 experiments, excluding 5 (respectively 4) unstable experiments for controllers (ii) and (iii). Controller~(i) is the best performing one since it enforces exact nonlinearity cancellation. Controller (iii) is shown to outperform controller (ii), although the same basis functions \eqref{eqn_ex_basisfunctions} were used, potentially suggesting that the region of attraction of controller (iii) is larger compared to (ii). This can be attributed to the fact that larger levels of PE were attained using multiple experiments.
	\section{Conclusion}\label{sec_conc}
We provided explicit formulas for inputs that guarantee PE for linear and classes of nonlinear systems. For Hammerstein and locally reachable nonlinear systems (including SISO flat systems), we showed how to guarantee collective PE of input- and/or state-dependent basis functions. These results are crucial for the application of recent data-driven control schemes that require such PE conditions to be satisfied.
\begin{table}[!t]
	\caption{Average level of PE and cumulative stabilization error.}
	\vspace{-1.5em}
	\begin{center}
		\begin{tabular}{|l|l|l|l|}
				\hline
				Compared value & Controller (i) & Controller (ii) & Controller (iii) \\ \hline
				avg $\sigma_{\textup{min}}(\mathcal{H}_1(\vartheta))$
				& 0.0083 & 0.0083 & 0.0399 \\ \hline
				$\sum_{i=0}^{T-1}\frac{1}{T}|x_{1,i}|$ 
				& 0.0469 & 0.0620 & 0.0553 \\ \hline
				$\sum_{i=0}^{T-1}\frac{1}{T}|x_{2,i}|$ 
				& 0.0231 & 0.0382 & 0.0314 \\ \hline
		\end{tabular}
		\vspace{-2.5em}
		\label{table_comparison}
	\end{center}
\end{table}

	\bibliographystyle{IEEEtran}
	\bibliography{references}

\begin{thebibliography}{10}
\providecommand{\url}[1]{#1}
\csname url@samestyle\endcsname
\providecommand{\newblock}{\relax}
\providecommand{\bibinfo}[2]{#2}
\providecommand{\BIBentrySTDinterwordspacing}{\spaceskip=0pt\relax}
\providecommand{\BIBentryALTinterwordstretchfactor}{4}
\providecommand{\BIBentryALTinterwordspacing}{\spaceskip=\fontdimen2\font plus
\BIBentryALTinterwordstretchfactor\fontdimen3\font minus
  \fontdimen4\font\relax}
\providecommand{\BIBforeignlanguage}[2]{{%
\expandafter\ifx\csname l@#1\endcsname\relax
\typeout{** WARNING: IEEEtran.bst: No hyphenation pattern has been}%
\typeout{** loaded for the language `#1'. Using the pattern for}%
\typeout{** the default language instead.}%
\else
\language=\csname l@#1\endcsname
\fi
#2}}
\providecommand{\BIBdecl}{\relax}
\BIBdecl

\bibitem{Ljung87}
L.~Ljung, \emph{System Identification: Theory for the User}, 2nd~ed.\hskip 1em
  plus 0.5em minus 0.4em\relax Prentice Hall, 1999.

\bibitem{astrom08}
K.~{\AA}str{\"o}m and B.~Wittenmark, \emph{Adaptive Control}, 2nd~ed.\hskip 1em
  plus 0.5em minus 0.4em\relax Dover Publications, 2008.

\bibitem{Green86}
M.~Green and J.~B. Moore, ``Persistence of excitation in linear systems,''
  \emph{Systems \& Control Letters}, vol.~7, no.~5, pp. 351--360, 1986.

\bibitem{Willems05}
J.~C. Willems, P.~Rapisarda, I.~Markovsky, and B.~L. {De Moor}, ``A note on
  persistency of excitation,'' \emph{Systems \& Control Letters}, vol.~54,
  no.~4, pp. 325--329, 2005.

\bibitem{Markovsky21}
I.~Markovsky and F.~Dörfler, ``Behavioral systems theory in data-driven
  analysis, signal processing, and control,'' \emph{Ann. Rev. in Control},
  2021.

\bibitem{Berberich20}
J.~Berberich and F.~Allgöwer, ``A trajectory-based framework for data-driven
  system analysis and control,'' in \emph{19th IEEE ECC}, 2020.

\bibitem{AlsaltiBerLopAll2021}
M.~Alsalti, J.~Berberich, V.~G. Lopez, F.~Allgöwer, and M.~A. Müller,
  ``Data-based system analysis and control of flat nonlinear systems,'' in
  \emph{60th IEEE CDC}, 2021, pp. 1484--1489.

\bibitem{Alsalti2022}
M.~Alsalti, V.~G. Lopez, J.~Berberich, F.~Allgöwer, and M.~A. Müller,
  ``Data-based control of feedback linearizable systems,'' \emph{IEEE
  Transactions on Automatic Control}, pp. 1--8, 2023.

\bibitem{DePersis22}
C.~De~Persis, M.~Rotulo, and P.~Tesi, ``Learning controllers from data via
  approximate nonlinearity cancellation,'' \emph{IEEE Transactions on Automatic
  Control}, pp. 1--16, 2023.

\bibitem{vanWaarde22}
H.~J. van Waarde, ``Beyond persistent excitation: Online experiment design for
  data-driven modeling and control,'' \emph{IEEE Control Systems Letters},
  vol.~6, pp. 319--324, 2022.

\bibitem{Yuan22}
Z.~Yuan and J.~Cortés, ``Data-driven optimal control of bilinear systems,''
  \emph{IEEE Control Systems Letters}, vol.~6, pp. 2479--2484, 2022.

\bibitem{DePersis21}
C.~{De Persis} and P.~Tesi, ``Designing experiments for data-driven control of
  nonlinear systems,'' \emph{IFAC-PapersOnLine}, 2021.

\bibitem{Alsalti2021c}
M.~Alsalti, V.~G. Lopez, J.~Berberich, F.~Allgöwer, and M.~A. Müller,
  ``Data-driven nonlinear predictive control for feedback linearizable
  systems,'' \emph{arXiv: 2211.06339}, accepted for 22nd IFAC WC, 2023.

\bibitem{vanWaarde20}
H.~J. van Waarde, C.~De~Persis, M.~K. Camlibel, and P.~Tesi, ``Willems’
  fundamental lemma for state-space systems and its extension to multiple
  datasets,'' \emph{IEEE Control Syst. Lett.}, vol.~4, no.~3, 2020.

\bibitem{Coulson22}
J.~Coulson, H.~J. van Waarde, J.~Lygeros, and F.~Dörfler, ``A quantitative
  notion of persistency of excitation and the robust fundamental lemma,''
  \emph{IEEE Control Systems Letters}, pp. 1--1, 2022.

\bibitem{Siami21}
M.~Siami, A.~Olshevsky, and A.~Jadbabaie, ``Deterministic and randomized
  actuator scheduling with guaranteed performance bounds,'' \emph{IEEE
  Transactions on Automatic Control}, vol.~66, no.~4, 2021.

\bibitem{XinJin19}
X.~Jin, R.~Bighamian, and J.-O. Hahn, ``Development and in silico evaluation of
  a model-based closed-loop fluid resuscitation control algorithm,'' \emph{IEEE
  Trans. Biomed. Eng.}, vol.~66, pp. 1905--1914, 2019.

\bibitem{HO66}
B.~HO and R.~E. K{\'a}lm{\'a}n, ``Effective construction of linear
  state-variable models from input/output functions,''
  \emph{at-Automatisierungstechnik}, vol.~14, no. 1-12, pp. 545--548, 1966.

\bibitem{Sun17}
C.~Sun and R.~Dai, ``Rank-constrained optimization and its applications,''
  \emph{Automatica}, vol.~82, pp. 128--136, 2017.

\bibitem{Markovsky13_SLRA}
I.~Markovsky and K.~Usevich, ``Structured low-rank approximation with missing
  data,'' \emph{SIAM Journal on Matrix Analysis and Applications}, vol.~34,
  no.~2, pp. 814--830, 2013.

\bibitem{Melkior20}
M.~Ornik, ``Guaranteed reachability for systems with unknown dynamics,'' in
  \emph{2020 59th IEEE CDC}, 2020, pp. 2756--2761.

\bibitem{Sontag13}
E.~D. Sontag, \emph{Mathematical Control Theory}, 2nd~ed., ser. Texts in
  Applied Mathematics.\hskip 1em plus 0.5em minus 0.4em\relax Springer New
  York, NY, 2013.

\bibitem{Christensen06}
O.~Christensen and K.~L. Christensen, ``Linear independence and series
  expansions in function spaces,'' \emph{The American Mathematical Monthly},
  vol. 113, no.~7, pp. 611--627, 2006.

\bibitem{MonacoNor1987}
S.~Monaco and D.~Normand-Cyrot, ``Minimum-phase nonlinear discrete-time systems
  and feedback stabilization,'' in \emph{26th IEEE CDC}, vol.~26, 1987, pp.
  979--986.

\end{thebibliography}

\end{document}